\RequirePackage{fix-cm}
\documentclass[smallextended]{svjour3}       % onecolumn (second format)
\smartqed  % flush right qed marks, e.g. at end of proof
\usepackage{graphicx}
\usepackage{amsmath}
\usepackage{amssymb}
\usepackage{mathtools}
\usepackage{booktabs}
\usepackage{tikz}
\newcommand{\R}{\mathbb{R}}
\newcommand{\Dt}{\delta t}

\begin{document}

\title{A Non-Intrusive Parallel-in-Time Adjoint Solver with the XBraid Library}
% \subtitle{Do you have a subtitle?\\ If so, write it here}

% \titlerunning{Adjoint Sensitivity Computation for XBraid}        % if too long for running head

\author{Stefanie G\"unther        \and
        Nicolas R. Gauger       \and
        Jacob B. Schroder
}

\authorrunning{S. G\"unther \and N.R. Gauger \and J.B. Schroder} % if too long for running head

\institute{S. G\"unther, N.R. Gauger \at
              TU Kaiserslautern\\
              Chair for Scientific Computing\\
               Paul-Ehrlich-Stra{\ss}e 34 -- 36\\
                67663 Kaiserslautern\\
                Tel.: +49-631-205-5637\\
              Fax: +49-631-205-3056\\
              \email{stefanie.guenther@scicomp.uni-kl.de}           %  \\
            % \emph{Present address:} of F. Author  %  if needed
           \and
           J.B. Schroder \at
           Center for Applied Scientific Computing\\
           Lawrence Livermore National Laboratory \\
           P.O. Box 808\\
           L-561, Livermore, CA 94551\\
           This work performed under the auspices of the U.S. Department of
           Energy by Lawrence Livermore National Laboratory under Contract
           DE-AC52--07NA27344, LLNL-JRNL-730159.  }

\date{Received: date / Accepted: date}

\maketitle

\begin{abstract}

  In this paper, an adjoint solver for the multigrid-in-time software library XBraid is presented. XBraid provides a non-intrusive approach for simulating unsteady dynamics on multiple processors while parallelizing not only in space but also in the time domain~\cite{xbraid-package}. It applies an iterative multigrid reduction in time algorithm to existing spatially parallel classical time propagators and computes the unsteady solution parallel in time. Techniques from Automatic Differentiation are used to develop a consistent discrete adjoint solver which provides sensitivity information of output quantities with respect to design parameter changes. The adjoint code runs backwards through the primal XBraid actions and accumulates gradient information parallel in time. It is highly non-intrusive as existing adjoint time propagators can easily be integrated through the adjoint interface. The adjoint code is validated on advection-dominated flow with periodic upstream boundary condition. It provides similar strong scaling results as the primal XBraid solver and offers great potential for speeding up the overall computational costs for sensitivity analysis using multiple processors.

\end{abstract}

\section{Introduction}\label{sec:intro}

With the rapid increase in computational capacities, computational fluid dynamics (CFD) has nowadays become a powerful tool to predict and analyze {steady and unsteady fluid flows. For unsteady dynamics, a numerical CFD simulation approximates the state $u(t)$ of a dynamical system (e.g. density, velocity or temperature), given problem specific design parameters denoted by $\rho$ (such as geometry, boundary and initial conditions or material coefficients),} while the state is determined through a set of unsteady partial differential equations (PDEs):
  \begin{align}
    \frac{\mathrm{d}u(t)}{\mathrm{d}t} &= g(u(t), \rho) \quad \forall  t \in (0,T) \\
    u(0) &= g^0.
\end{align}
{Here, the right hand side $g$ involves spatial derivatives, source terms, etc. as well as the design parameters $\rho$}. Accounting for the unidirectional flow of information in the time domain, an approximation to the unsteady PDE solution is typically evolved forward in time in a step-by-step manner applying a time marching algorithm~\cite{blazek2005computational,Peric:2002,tucker2014unsteady}.
Starting from the initial condition, these schemes march forward in discrete time steps
% with
%   \begin{align}\label{timemarchingscheme}
%     \vec u^{i} = \Phi^i(\vec u^{i-1}, \rho), \quad \text{for }i=1,\ldots, N.
%   \end{align}
% An application of the time stepper $\Phi^i$ usually involves
% these schemes compute approximations to the state at discrete time steps,
applying nonlinear {iterations in space} to approximate a pseudo-steady state at each time step as for example in the dual-time stepping approach~\cite{jameson1991dualtime}.

However, in many applications, the primal flow is not the only computation of particular interest. The ability to compute sensitivities is also needed to determine the influence of design changes to some objective function $J$ that {computes} the time-average of instantaneous quantities of the flow dynamics
\begin{align}
  J(u, \rho) = \frac{1}{T}\int_0^T f(u(t), \rho) \mathrm{d}t.
\end{align}
The ability to compute sensitivities can improve and enhance the simulation tool, as for example through parameter estimation for validation and verification purposes~\cite{navon1998practical}, error estimation and adaptive grid refinement~\cite{venditti2002grid,pierce2000adjoint,giles2002adjoint} or uncertainty quantification techniques~\cite{beyer2007robust}. Further, it broadens the application range from pure simulation to optimization as for example in an optimal control or shape optimization framework~\cite{jameson1988aerodynamicdesign,nadarajah2007optimum,mohammadi2010applied}.

If the number of independent design parameters is large, the adjoint approach for sensitivity computations is usually preferred since its computational cost does not scale with the design space dimension~\cite{pironneau_1974,giles1997adjoint}.
In that approach, only one additional adjoint equation needs to be solved in order to set up the gradient of the objective function $J$ with respect to $\rho$.
% In that approach, perturbations of the output objective function are traced back to perturbation of the independent input parameters in order to set up the gradient efficiently as will be presented in Section \ref{sec:AdjointSensivity}. F
For unsteady time marching schemes, solving the adjoint equation involves a reverse time integration loop that propagates sensitivities backwards through the time domain
% \begin{align}\label{adjointtimemarchingscheme}
%   \bar{\vec u}^{i} &= \nabla_{\vec u^{i}}f(\vec u^{i}, \rho) + \left(\partial_{\vec u^{i}} \Phi^{i+1}(\vec u^{i}, \rho)\right)^T\bar{\vec u}^{i+1} \quad \text{for } i=N, \dots, 1
% \end{align}
starting from a terminal condition~\cite{nadarajah2007optimum,economon2013unsteady,rumpfkeil2007general,mavriplis2008solution,nielsen2010discrete}. Hence evaluating the gradient is a rather computationally expensive task as it involves a forward loop over the time domain to approximate the PDE solution followed by a backwards time marching loop for the adjoint.

One way to mitigate this costly procedure of forward and backward sweeps over
the time domain is to parallelize the time dimension.  Overall, the need for
parallel-in-time methods is being driven by changes in computer architectures,
where future speedups will be available through greater concurrency, not faster
clock speeds, which are stagnant.  Previously, ever faster clock speeds made it
possible to speedup fixed-size sequential time stepping problems
{(strong scaling),}
and also to refine problems in time and space without increasing the
wall-clock time {(weak scaling)}.  However with stagnate clock speeds, speedups
will only be possible through the increasing concurrency of future
architectures.  Thus, to overcome this serial time stepping bottleneck, time
stepping codes must look for new parallelism from the time dimension.  This
bottleneck is particularly acute for adjoint solvers because of the many
sweeps over the time domain.

Research on parallel-in-time methods goes at least back to the 1964 work by
Nievergelt~\cite{Nievergelt_1964}.  Since then, a variety of approaches have
been developed, including parareal, which is perhaps
the most popular parallel-in-time method~\cite{parareal-2001}. For a recent literature review and a
gentle introduction, see~\cite{Ga2015}.  The method used here is multigrid
reduction in time (MGRIT) and is based on multigrid
reduction~\cite{mgrit-2014,ries_trottenberg_1979_mgr,ries_etal_1983_mgr}.  MGRIT is
relatively non-intrusive on existing codes and when restricted to two grid
levels, is equivalent to parareal~\cite{mgrit-2014,GV_2007_parareal}. This
multilevel distinction is important, because it allows for optimal parallel
communication behavior, as opposed to a two-level scheme, where the size of the
coarse-level limits concurrency.  This work uses XBraid,
an open source implementation of MGRIT developed at Lawrence Livermore National
Lab~\cite{xbraid-package}.

Previous time parallel approaches to optimization include multiple shooting
methods~\cite{He2005}, parareal approaches applied to a reduced Hessian~\cite{DuSaSch2005,mathew2010analysis}, and the Schwarz preconditioner approaches in~\cite{BaSt2015,KwGa2016,kwok2017time}. {Further, \cite{GoetschelMinion2017} uses the time-parallel PFASST approach to solve the state and adjoint equations for parabolic optimal control problems.} While these previous approaches were effective, applying XBraid
for optimization offers some possible advantages, e.g., XBraid offers multilevel scalability,
and the potential for greater parallelism over the Schwarz approaches.  The
parallelism for the Schwarz approaches is limited by the size of fine-grid ``subdomains''
which must each be solved sequentially, but for XBraid, the size
of such subdomains is determined by the temporal coarsening factor, which can be
as small as 2.

%
%An additional reason that optimization is a field ripe for use with
%parallel-in-time is memory usage.  A typical concern of parallel-in-time is its
%increased memory usage; however, the classical optimization approach often
%checkpoints data on a regular basis during the forward sweep to save
%re-computation during the backward sweep. In short, the extra storage required
%for parallel-in-time is often already being done.

When considering the intersection of optimization, CFD, and
parallel-in-time, a
prerequisite is that parallel-in-time be effective for standard
fluid dynamics problems, and recent advances indicate that parallel-in-time can
indeed be an effective tool here.
In \cite{FarCha2003}, the Parallel Implicit Time-Integration Algorithm (PITA)
was proposed, which enjoys success for highly advective problems, targeting in
part fluid dynamics~\cite{RuKr2012}.  It essentially uses a parareal framework,
but creates a massive Krylov space to stabilize and accelerate the method.
While effective, the memory usage of this algorithm is discouraging.  When
considering unmodified parareal, the works \cite{FiHeMa2005,SteRuSpKr2015,RuKr2014} show that parareal can be effective
for fluids problems when sufficient physical (not artificial) diffusion is
present.  Similar success with XBraid was shown for vortex
shedding~\cite{FalKatz2014}.  Lastly, the work \cite{Do2016} showed a path forward for XBraid and
purely advective problems that relies only on artificial
dissipation.

The XBraid package will be presented in Section
\ref{sec:XBraid}. Section \ref{sec:AdjointSensivity} reviews the adjoint approach
for computing sensitivities. Then, in Section \ref{sec:XBraid_adjoint}, an approach to
modifying XBraid in order to integrate existing adjoint time stepping codes
into the XBraid framework is developed. Exploiting techniques from Automatic
Differentiation (AD)~\cite{griewank2008evaluating}, the primal XBraid
iterations are enhanced by an adjoint iteration that runs backwards through the
XBraid actions and computes consistent discrete adjoint sensitivities. Similar
to the primal XBraid solver, the resulting adjoint solver is non-intrusive in
the sense that existing adjoint simulation codes can easily be integrated
through the extended user interface. Finally, the adjoint code is validated in Section \ref{sec:Numerics} by applying it to advection-dominated flow with a periodic upstream boundary.

\section{Multigrid Reduction in Time using XBraid} \label{sec:XBraid}

To describe the MGRIT algorithm implemented by XBraid, let $u(t)$ be
the solution to a time-dependent problem on the
interval $[0,T]$.
Let $t_i = i \Dt, i = 0, 1, ..., N$ be a time grid of that interval with
spacing $\Dt = T/N$, and let $\vec u^{i}\in \R^n$ be an approximation of $u(t_i)$.

Classical time marching schemes successively compute $\vec u^i$ based on information at previous time steps $\vec u^{i-1}, \vec u^{i-2}, \dots,$ and the design $\rho$ where {$\rho\in\R^p$ in a discretized setting}. Implicit methods are often prefered due to better stability properties leading to nonlinear equations at each time step that need to be solved iteratively for $\vec u^i$. Wrapping these iterations into a time stepper {$\Phi^i\colon \R^n\times \R^p \to \R^n$}, a general one-step time marching scheme can be written as
\begin{eqnarray} \label{eqn-st-system}
   \vec u^{i} = \Phi^{i}(\vec u^{i-1}, \rho), \quad \mbox{for } i = 1,2, ..., N,
\end{eqnarray}
with given initial condition {$\vec u^0 = g^0$}. Application to multi-step methods which involve more than one previous time step can be found in \cite{falgout2015multigrid} and are based on a reformulation into a block one-step scheme.

{We consider first the linear case (only as motivation), so
that $\Phi^{i}(\vec u^{i-1}, \rho)$ is equivalent to the matrix-vector product
$\Phi^i_{\rho} \vec u^{i-1}$, yielding $\vec u^i = \Phi^i_{\rho} \vec u^{i-1} + \vec
g^i_{\rho}$.  Thus, time stepping is equivalent to a forward sequential solve
of the block lower bidiagonal system $A \vec u = \vec g$,\footnote{In this
work, the application of $\Phi^i_{\rho}$ represents the approximate inversion
of an operator (implicit scheme), but in principle, explicit schemes may be
considered as well.} }
\begin{equation}\label{eqn-phi-system}
   \left(
   \begin{array}{cccc}
        I       &        &         & \\
          -\Phi^{1}_{\rho} & I      &         & \\
                    & \ddots & \ddots  & \\
                              &        & -\Phi^{N}_{\rho} & I
   \end{array}
   \right)
   \left(
   \begin{array}{c}
    \vec u^{0} \\
    \vec u^{1} \\
    \vdots \\
    \vec u^{N}
   \end{array}
   \right)
   =
   \left(
   \begin{array}{c}
     g^{0} \\
     g^{1}_{\rho} \\
    \vdots \\
     g^{N}_{\rho}
   \end{array}
   \right).
\end{equation}
The MGRIT algorithm replaces the $O(N)$ sequential time stepping algorithm
with a highly parallel $O(N)$ multigrid algorithm~\cite{oosterlee_book}.
%The MGRIT method is derived from approximate block cyclic reduction methods,
%which are well-known methods for tridiagonal systems.
Essentially, a sequence of coarser temporal grids are used to accelerate
the solution of the fine grid (\ref{eqn-phi-system}).  The MGRIT algorithm
easily extends to nonlinear problems with the Full Approximation Storage (FAS)
nonlinear multigrid method~\cite{Brandt_1977},
{which is the same nonlinear multigrid cycling used by PFASST \cite{emmett2012}.}
MGRIT solves the discrete space-time system
in equation (\ref{eqn-st-system}); i.e., MGRIT converges to the same solution
produced by sequential time stepping on the finest grid.  Lastly, while
sequential time stepping and MGRIT are optimal, the constant for MGRIT is
higher,
{ often by a factor of 10 or 20 for a straight-forward application of MGRIT,}
(an optimal direct sequential method has been replaced with an
optimal iterative method).  This creates a crossover point where a certain
number of resources are required to overcome the extra computational work of MGRIT.

\begin{figure}
\center
\usetikzlibrary{decorations.pathreplacing}
\begin{tikzpicture}
\path[draw] (0,0) -- (10,0);

% Draw Fine Points
\foreach \x in {0,1,...,3}
   \path[draw, line width = 1pt] (\x/2,0.1) -- (\x/2,-0.1) node [below] {$t_{\x}$};
\path[draw, line width = 1pt] (2,0.1) -- (2,-0.1) node [below] {$...$};
\path[draw, line width = 1pt] (2.5,0.1) -- (2.5,-0.1) node [below] {$t_m$};
\foreach \x in {4,5,...,20}
   \path[draw, line width = 1pt] (\x/2,0.1) -- (\x/2,-0.1);

% Draw Coarse Points
\path[draw, line width = 2pt, color = red] (0,-0.2) -- (0,0.2) node [above] {$T_0$};
\path[draw, line width = 2pt, color = red] (2.5,-0.2) -- (2.5,0.2) node [above] {$T_1$};
\path[draw, line width = 2pt, color = red] (5,-0.2) -- (5,0.2) ;
\path[draw, line width = 2pt, color = red] (7.5,-0.2) -- (7.5,0.2) ;
\path[draw, line width = 2pt, color = red] (10,0.2) -- (10,-0.2) node [below] {$\textcolor{black}{t_{N_t}}$};
\draw [decoration={brace}, decorate, line width = 1pt] (5,0.23) -- (7.5,0.23) node [above, pos=0.5] {$\Delta T = m \delta t$};
\draw [decoration={brace,mirror}, decorate, line width = 1pt] (5.5,-0.22) -- (6,-0.22) node [below, pos=0.5] {$\delta t$};
\end{tikzpicture}
\caption{Fine grid ($t_i$) and coarse grid ($T_j$) for coarsening factor
$m=5$.  The coarse grid induces a decomposition of the fine grid into
$C$-points (red) and $F$-points (black).}
\label{fig-grid}
\end{figure}
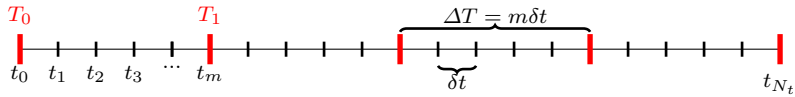

The MGRIT algorithm forms its sequence of coarse time grids from the original fine grid by
successively coarsening with factor $m>1$.  When a time grid is coarsened, it is
decomposed into two sets called $C$-points (points that will go on to form the
next coarser grid) and $F$-points (points that exist only on the fine grid).
Figure~\ref{fig-grid} illustrates an example of this.  This decomposition then
in turn, induces the relaxation method and coarse-grid correction step, which
together form the basis of a multigrid method.

Relaxation is a block Jacobi method that alternates between the $F$- and
$C$-points.  An F-relaxation updates $F$-point values $\vec u^{i}$ on interval
$(T_j, T_{j+1})$ by propagating the $C$-point value $\vec u^{mj}$ with $\Phi^{i}$.
Each $F$-interval is independent, and can be computed in parallel.  Likewise,
C-relaxation uses $\Phi^{mj}$ to update each $C$-point value $\vec u^{mj}$ based
on the previous $F$-point value.  These updates can also be carried out in
parallel.  Thus, FCF-relaxation corresponds to three relaxation sweeps over
$F$-points, $C$-points and $F$-points.

The coarse time grid problem is constructed by rediscretizing the problem using
only the red $C$-points, as depicted in Figure \ref{fig-grid}.  This
corresponds to eliminating the $F$-points in equation (\ref{eqn-phi-system}),
followed by the substitution of a cheaper coarse time step operator for the
exact fine-grid time step operator.  Injection is used to transfer vectors
between grids (see \cite{mgrit-2014} for more details).  When the algorithm is
done recursively, using FAS, a multilevel algorithm capable of solving complex
nonlinear problems is produced (see \cite{FalKatz2014,FaMaOnSc2017}).  Standard
multigrid cycling strategies, e.g., V-cycles and F-cycles (see
\cite{oosterlee_book}), can then be applied.
% --> Could insert picture of V- and F-cycles, as well as an actual algorithm,
% but I think space argues against it.

The FAS nonlinear multigrid method is a general nonlinear solver method, and as
such, the solution of the global space-time system is a fixed-point of the MGRIT algorithm.
For instance, the paper \cite{Do2016} shows that the error propagator of
MGRIT is a contraction for many linear cases.  Thus, letting $\vec u_k = \left(\vec u^1_k, \dots, \vec u^N_k\right) \in \R^{Nn}$ denote the
solution to the space-time system after $k$ MGRIT iterations, and the
operator $H\colon \R^{Nn} \times \R^p\to \R^{Nn}$ represent the application of
one MGRIT iteration, the overall solution process can be represented as
\begin{align}
   \label{eqn:fixed_point}
   \vec u_{k+1} = H(\vec u_k, \rho),
\end{align}
where $\vec u_{k+1}$ is converging to a fixed-point $\vec u = H(\vec u, \rho)$ with $\vec u^i = \Phi^i(\vec u^{i-1}, \rho)$ at all time steps $i=1, \dots, N$.

%In subsequent sections, the superscript $^{i}$ notation which denotes
%time step number is dropped for ease of reading.

\subsection{XBraid Interface}\label{subsec:XBraidInterface}

The XBraid interface places a high value on
non-intrusiveness and is designed to wrap existing time stepping codes.  The
goal is for the user to generate only a small amount of new wrapper code that
then yields a new parallel-in-time capability.  XBraid is written in C, but
also has an object oriented C++ interface and a Fortran90 interface.
The time parallelism is handled with MPI.

To use XBraid, the user is responsible for defining and implementing the
following.  See \cite{xbraid-package} and the associated User's Manual for more
details.
\begin{itemize}

   \item A user $\mathtt{App}$ (application) structure must be defined and is globally
      available to the user in all wrapper routines such as
      $\mathtt{my\_Step}$.  This C-style structure generally holds time
      \emph{independent} information such as the time-grid definition, MPI
      communicators as well as the design variable $\rho$ and the current time-averaged objective function $J$.

   \item A state $\mathtt{vector}$ structure must be defined so that it captures
      a complete solution
      snapshot at a single time point.  This C-style structure generally
      holds time \emph{dependent} information, e.g., unknown values $\vec u^i$ and
      time-dependent spatial mesh information.

   \item The key user-written routine is the $\mathtt{my\_Step}$ function which takes as
      input a state $\mathtt{vector}$ at time step $i-1$ and advances it to
      time step $i$.  This function is called everywhere, i.e., on coarse and fine time-grids
      for large and small time step sizes.  This function wraps the user's
      existing time stepping routine and defines $\Phi$, i.e.
      $\vec u^{i} = \Phi^i(\vec u^{i-1},\rho).$

   \item The user also writes the output routine $\mathtt{my\_Access}$.
      This function is called after
      each XBraid iteration to access and possibly modify the state vectors and to pass information on the instantaneous quantity of interest $f(\vec u^i, \rho)$ to the user $\mathtt{App}$.  It is called individually for each time step value
      on the finest time grid.

   \item The user must further define a variety of wrapper functions that provide XBraid with additional information for handling the $\mathtt{vector}$ structure. For instance, $\mathtt{my\_Init}$ and
      $\mathtt{my\_Clone}$ must initialize and clone a $\mathtt{vector}$,
      respectively.  The wrapper $\mathtt{my\_Sum}$ must compute a vector sum
      of two $\mathtt{vectors}$, i.e., $\vec v^i = \alpha \vec u^i + \beta \vec v^i$.
      Other functions must free a $\mathtt{vector}$ structure and take
      the norm of a $\mathtt{vector}$, typically an Euclidean norm.
   \item Additionally, the user must define how to pack and unpack an MPI buffer
      so that $\mathtt{vectors}$ can be sent between processors.  These
      $\mathtt{my\_BufPack}$ and $\mathtt{my\_BufUnpack}$ functions are
      generally straightforward and require the user to flatten a
      $\mathtt{vector}$ into a buffer, and then unpack that buffer into a
      $\mathtt{vector}$.

      %$\mathtt{my\_Init}$, $\mathtt{my\_Clone}$, $\mathtt{my\_Free}$,
      %$\mathtt{my\_Sum}$, $\mathtt{my\_SpatialNorm}$, $\mathtt{my\_BufPack}$,
      %$\mathtt{my\_BufUnpack}$, $\mathtt{my\_BufSize}$

\end{itemize}

\section{Adjoint Sensitivity Computation}\label{sec:AdjointSensivity}
The presented XBraid solver determines a global space-time solution by iteratively solving the fixed-point equation
\begin{align}\label{adj:fixedpointequation}
 \vec u = H(\vec u, \rho),
\end{align}
for a fixed design $\rho$. Small changes in the design variable $\delta \rho$ result in a perturbed state $\delta \vec u$ which induces a perturbation in the objective function $\delta J$ both explicitly from the design itself and implicitly through the solution of fixed-point equation $\eqref{adj:fixedpointequation}$. This change is quantified by the sensitivity (or derivative) of $J$ with respect to $\rho$.

One way to determine the sensitivity numerically is the finite difference approach which approximates the derivative for design perturbations with unity directions $e_j\in \R^p$ by evaluating
\begin{align}
  \frac{\mathrm{d}J}{\mathrm{d}\rho_j}\ \approx \frac{J(\vec u_{\rho+\epsilon e_j}, \rho+\epsilon e_j) - J(\vec u, \rho)}{\epsilon} \quad \forall \, j=1, \dots, p,
\end{align}
for small $\epsilon>0$, where $\vec u_{\rho+\epsilon e_j}$ is the space-time solution of the fixed-point equation for the perturbed design parameter. However, a full unsteady simulation would be necessary to compute $\vec u_{\rho+\epsilon e_j}$ such that the computational cost for computing the full gradient of $J$, which consists of the derivatives in all unity directions in $\R^p$, grows proportional to the design space dimension $p$.

\subsection{Adjoint sensitivities via Lagrangian formalism}

A well-established alternative approach to computing the gradient of $J$ with respect to $\rho$ is the adjoint method~\cite{lions1971optimal,pironneau_1974,Giles2000}. In that approach, an additional equation is derived by differentiating $H$ and $J$ partially with respect to $\vec u$ and $\rho$. Its solution, the so-called adjoint variable, is then used to determine the gradient of $J$ for a computational cost that does not scale with the number of design parameters. The adjoint equation can be derived using an augmented function (the so-called Lagrangian function)
\begin{align}
  L(\vec u, \bar{\vec u}, \rho) \coloneqq J(\vec u, \rho) + \bar{\vec u}^T \left(H(\vec u, \rho) - \vec u\right),
\end{align}
which adds a multiple of the fixed-point equation \eqref{adj:fixedpointequation} to the objective function with a multiplier $\bar{\vec u} \in \R^{Nn}$ . If $\vec u$ is a solution to the fixed-point equation, the derivatives of $J$ and $L$ with respect to $\rho$ coincide and are computed from the chainrule as
 \begin{align}
   \frac{\mathrm{d}J}{\mathrm{d}\rho} = \frac{\mathrm{d}L}{\mathrm{d}\rho}
        &= \frac{\partial J}{\partial \vec u}\frac{\mathrm{d}\vec u}{\mathrm{d}\rho} + \frac{\partial J}{\partial \rho} + \bar{\vec u}^T\left( \frac{\partial H}{\partial \vec u}\frac{\mathrm{d}\vec u}{\mathrm{d}\rho} + \frac{\partial H}{\partial \rho} - \frac{\mathrm{d}\vec u}{\mathrm{d}\rho} \right)\\
        &= \frac{\partial J}{\partial \rho} + \bar{\vec u}^T\frac{\partial H}{\partial \rho} + \left( \frac{\partial J}{\partial \vec u} + \bar{\vec u}^T\frac{\partial H}{\partial \vec u} -\bar{\vec u}^T \right)\frac{\mathrm{d}\vec u}{\mathrm{d}\rho}. \label{adj:chainrule}
 \end{align}
 The term $\mathrm{d}\vec u / \mathrm{d} \rho$ represents the sensitivity of the flow solution with respect to design changes whose numerical computation is extremely expensive as it would require $p$ primal flow simulations in a finite difference setting. The adjoint approach avoids these computations by choosing the multiplier $\bar{\vec u}$ in such a way that terms containing these sensitivities add up to zero: Choosing $\bar{\vec u}$ such that it satisfies the so-called adjoint equation
 \begin{align}\label{adj:adjointequation}
   \bar{\vec u} = \nabla_{\vec u}J + \left( \frac{\partial H}{\partial \vec u}\right)^T \bar{\vec u},
 \end{align}
the gradient of $J$ is then given by
\begin{align}\label{adj:reducedgradient}
  \nabla J &= \left(\frac{\mathrm{d}J}{\mathrm{d}\rho}\right)^T = \nabla_{\rho} J + \left(\frac{\partial H}{\partial \rho}\right)^T \bar{\vec u},
\end{align}
according to equation $\eqref{adj:chainrule}$, {where $\nabla$ denotes transposed derivative vectors.}
Instead of solving the fixed-point equation $p$ times for each unity direction in the design space, the adjoint approach requires only the solution of one additional (adjoint) equation for the adjoint variable $\bar{\vec u}$ in order to determine the gradient of $J$ with respect to $\rho$.

\subsection{Simultaneous piggyback iteration for the primal and adjoint variable}
The adjoint equation is linear in the adjoint variable. Hence, if $H$ is contractive,
% with
% \begin{align}
%   \left\|\frac{\partial H}{\partial \vec u} \right\| \leq \rho < 1
% \end{align}
% for all points of interest,
it can be solved with the linear fixed-point iteration
\begin{align}
  \bar{\vec u}_{k+1} = \nabla_{\vec u}J(\vec u, \rho) + \left( \frac{\partial H(\vec u, \rho)}{\partial \vec u}\right)^T \bar{\vec u}_k \quad \text{for } k=0,1,\dots,
\end{align}
starting from $\bar{\vec u}_0 = \vec 0$ where the right hand side is evaluated at a feasible state $\vec u$ that satisfies \eqref{adj:fixedpointequation}. However, it is shown in \cite{griewank2002reduced} that the adjoint iteration can also be performed simultaneously with the primal iterations in the following \textit{piggyback} approach:
\begin{align}
  \text{For } k=0,1,\dots&: \notag \\
  \vec u_{k+1} &=  H(\vec u_k, \rho) \label{piggyback_itertion_primal} \\
   \bar{\vec u}_{k+1}  &= \nabla_{\vec u}J(\vec u_k, \rho) + \left( \frac{\partial H(\vec u_k, \rho)}{\partial \vec u}\right)^T \bar{\vec u}_k\label{piggyback_itertion_adjoint}.
\end{align}

This approach is especially attractive for simultaneous optimization algorithms as for example the One-shot method where small changes to the design are introduced in each piggyback iteration based on the current evaluation of the gradient~\cite{griewank2006projected,bosse2014,gauger2012automated,gunther2016simultaneous}.
% The One-shot method has gained popularity in optimization with steady-state PDEs where it offers speedup over classical optimization methods and is under active development for unsteady PDE-constrained optimization.
For the case considered here of a fixed design $\rho$, the piggyback iterations converge to the feasible primal $\vec u$ and adjoint solution $\bar{\vec u}$ at the same asymptotic convergence rate determined by the contraction rate of $H$. Due to the dependency of the adjoint solution on a feasible state, the adjoint variable is expected to lag a bit behind the primal iterates which has been analyzed in \cite{griewank2004timelag}. If the contractivity assumption is slightly violated such that eigenvalues close to or even outside the unit sphere are present, the adjoint iteration can be stabilized using the recursive projection method (RPM), see for example \cite{shroff1993stabilization,albring2017rpm}.

\section{The adjoint XBraid solver}\label{sec:XBraid_adjoint}

While the primal XBraid solver provides an iteration for updating the primal state as in \eqref{piggyback_itertion_primal}, an adjoint XBraid solver has been developed that updates the adjoint variable as in \eqref{piggyback_itertion_adjoint} and evaluates the gradient at the current iterates. For current primal and adjoint input variables $\vec u_k, \bar{\vec u}_k$, it returns
\begin{align} \label{AD:adjoint_and_reducedgradient}
  \begin{bmatrix}\bar{\vec u}_{k+1} \\ \bar \rho \end{bmatrix}=
    \begin{bmatrix} \nabla_{\vec u} J(\vec u_k, \rho) + \partial_{\vec u}H(\vec u_k, \rho)^T\bar{\vec u}_k  \\
    \nabla_{\rho} J(\vec u_k, \rho) + \partial_{\rho} H(\vec u_k, \rho)^T\bar{\vec u}_k
  \end{bmatrix},
\end{align}
where $\bar{\vec u}_{k+1}$ is the new adjoint iterate and $\bar{\rho}$ holds the current approximation of the sensitivity of $J$ with respect to $\rho$ as in \eqref{adj:reducedgradient}.
% \subsection{Differentiation of XBraid}
Since $H$ and $J$ refer to numerical algorithms, their partial derivatives are computed by adopting techniques from Automatic Differentiation (AD)~\cite{griewank2008evaluating,naumann2012art}.

AD is a set of techniques that modify the semantics of a computer program in such a way that it not only computes the primal output but also provides sensitivities of the outputs with respect to input variations. It relies on the fact that any computer program for evaluating a numerical function can at runtime be regarded as a concatenation of elemental operations whose derivatives are known.  Thus, the derivative of the output can be computed by applying the chain rule to the elemental operations.
Two modes are generally distinguished: while the forward mode computes directional derivatives, the reverse mode returns transposed matrix-vector products of the sensitivities which is hence the method of choice for setting up the adjoint iteration.
% In the reverse mode of AD (often also called {adjoint mode}), the chain rule is evaluated from left to the right propagating sensitivities backwards through the computational graph.

To be more precise, consider a numerical algorithm for evaluating $z=F(x)$ with input $x\in \R^n$ and output $z\in \R^m$ which at runtime is a concatenation of elemental operations denoted by $h_l$ such that
\begin{align}
  z = F(x) = (h_L \circ h_{L-1} \circ \dots \circ h_1)(x).
\end{align}
 Given a vector $\bar z\in \R^m$, the reverse mode computes the sensitivity $\bar x \in \R^n$ with
  \begin{align}\label{AD:reversesensitivity}
    \bar x  = \left(\frac{\partial F(x)}{\partial x}\right)^T \bar z,
  \end{align}
  by applying the chain rule to the elemental operations. First the primal output itself is evaluated in a forward loop
  \begin{align} \label{AD:forwardloop}
    v_{l+1} = h_l(v_l) \quad l=0, \dots, L-1,
  \end{align}
  with $v_0 = x$, which yields $z=v_L$. Then a reverse loop evaluates and concatenates the local sensitivity of the elemental operations
  \begin{align}\label{AD:reverseloop}
    \bar v_l = \left( \frac{\partial h_l(v_l)}{\partial v_l} \right)^T \bar v_{l+1} \quad l=L-1,\dots, 0,
  \end{align}
using the reverse input $\bar v_{L} = \bar z$, which yields the transposed sensitivity product $\bar x = \bar v_0$ as in \eqref{AD:reversesensitivity}.

Following the above methodology, the adjoint XBraid solver is constructed by concatenating local sensitivities of elemental operations that produce the output $z=(\vec u, J) = (H(\vec u, \rho), J(\vec u, \rho)) = F(x)$ from the input $x=(\vec u, \rho)$ which yields
\begin{align}
  \bar x =
  \begin{pmatrix} \bar{\vec u} \\ \bar \rho
  \end{pmatrix}
   = \begin{pmatrix}
     \frac{\partial H}{\partial \vec u} & \frac{\partial H}{\partial \rho} \\ \frac{\partial J}{\partial \vec u} & \frac{\partial J}{\partial \rho}
   \end{pmatrix}^T
   \begin{pmatrix}
     \bar{\vec u} \\ \bar J
   \end{pmatrix} = \left(\frac{\partial F}{\partial x}\right)^T \bar z.
\end{align}
Choosing the reverse input vector $\bar z = (\bar{\vec u}, \bar J)= (\bar{\vec u}_k, 1)$, this produces the desired sensitivities as in \eqref{AD:adjoint_and_reducedgradient}.  Since XBraid modifies its input solely by calling the user-defined interface routines as introduced in Section \ref{subsec:XBraidInterface}, these interface routines are identified with the elemental operations $h_l$. The primal XBraid iteration manages the control flow of these actions and performs a forward loop for computing the primal output as in \eqref{AD:forwardloop}. The adjoint XBraid solver marches backwards through the same control flow in a reverse loop as in \eqref{AD:reverseloop} calling the differentiated routines $\bar h_l$ that evaluate the local sensitivities. A Last-In-First-Out (LIFO) data structure, called \textit{action tape} therefore memorizes the control flow of the primal actions $h_l$. The actions are overloaded in such a way, that they still calculate their primal output $v_{l+1}$ but as a side effect, they record themselves, their arguments $v_l$ and a pointer to the corresponding intermediate adjoint $\bar v_{l+1}$ on the action tape.
After the execution of a primal XBraid iteration, the reverse adjoint loop pops the elements from the action tape and calls the corresponding differentiated action $\bar h_l$ which updates the corresponding intermediate adjoint and the gradient with the correct local sensitivity information.
This essentially creates a separate adjoint XBraid solver that reverses through the primal XBraid actions collecting sensitivities. Thus, an equation of the form \eqref{eqn-phi-system} is solved but going backwards in time.

According to the reverse mode of AD, the local sensitivities are transposed matrix-vector products of the sensitivities of $h_l$ multiplied with the intermediate adjoint vectors. Table \ref{tab:diffXBraidActions} lists the original interface routines and their differentiated versions. Note, that the differentiated routines perform updates of the intermediate variables using the ``$a\mathrel{+}=b$'' notation for assigning ``$a \leftarrow a+b$'' rather then overwriting ``$a \leftarrow b$''.
% This results from the fact that only local input variables rather than the complete program state $v_l$ are stored on the tape in the current implementation.
This results from the fact that instead of storing the entire program state $v_l$ during the forward loop \eqref{AD:forwardloop}, only the local input variables for each action are pushed to the tape. Hence, it is possible that their intermediate adjoint variables are modified elsewhere in the code which necessitates updates rather than overwriting the variables (see e.g. \cite{griewank2008evaluating} for implementational details on the reverse mode of AD).

% from an improved implementation of the reverse mode that avoids storing the complete program state $v_l$ during the forward loop \eqref{AD:forwardloop} but rather stores only the local input variables (see e.g. \cite{griewank2008evaluating} for implementational details on reverse mode AD).

While the primal memory de-/allocation is managed by calling the primal user routines $\mathtt{my\_Init}$ and $\mathtt{my\_Free}$, memory management of the intermediate adjoint variables is automated by the use of shared pointers with reference counting. The current implementation uses the smart pointer class $\mathtt{std::shared\_ptr}$ defined in the C++11 standard library~\cite{josuttis2012c}. It counts the number of pointer copies to a specific object and destroys it automatically as soon as the last reference is removed.

% \begin{table}
% % table caption is above the table
% \caption{Table example}\label{tab:1}
% % For LaTeX tables use
% \begin{tabular}{lll}
% \hline\noalign{\smallskip}
% first & second & third  \\
% \noalign{\smallskip}\hline\noalign{\smallskip}
% number & number & number \\
% number & number & number \\
% \noalign{\smallskip}\hline
% \end{tabular}
% \end{table}

  \begin{table}
    \caption{Original and differentiated XBraid actions.}
    \label{tab:diffXBraidActions}
     \def\arraystretch{1.2}
    \begin{tabular}{lll}
      \hline\noalign{\smallskip}
      \textbf{Action} & \textbf{Original $h_l$} & \textbf{Differentiated $\bar h_l$} \\
      \noalign{\smallskip}\hline\noalign{\smallskip}
      % \midrule
      $\mathtt{my\_Step}$ & ${\vec u}^{i} = \Phi^i({\vec u}^{i-1}, \rho)$  &
      $\bar {\vec u}^{i-1} \mathrel{+}= (\partial_{{\vec u}^{i-1}} \Phi^i({\vec u}^{i-1},\rho))^T \bar {\vec u}^{i}$ \\
      & & $\bar \rho \mathrel{+}= (\partial_{\rho} \Phi^i({\vec u}^{i-1},\rho))^T \bar {\vec u}^{i}$ \\

      $\mathtt{my\_Access}$ & $f({\vec u}^i, \rho)$  & $\bar {\vec u}^i \mathrel{+}= \nabla_{{\vec u}^i} f({\vec u}^i,\rho)$ \\
      & & $\bar \rho \mathrel{+}= \nabla_{\rho} f({\vec u}^i,\rho)$ \\

      $\mathtt{my\_Clone}$ & ${\vec v}^i = {\vec u}^i$ & $\bar {\vec u}^i \mathrel{+}= \bar{\vec v}^i$ \\
       & & $\bar{\vec v}^i = 0$ \\

      $\mathtt{my\_Sum}$ & $\vec v^i = \alpha {\vec u}^i+ \beta \vec v^i$ & $\bar {\vec u}^i \mathrel{+}= \alpha \bar{\vec v}^i$ \\
      & & $\bar {\vec v}^i = \beta \bar {\vec v}^i$ \\

      $\mathtt{my\_BufPack}$ & MPI\_Send(${\vec u}^i$) & MPI\_Recv($\bar {\vec u}^i$) \\
      $\mathtt{my\_BufUnPack}$ & MPI\_Recv(${\vec u}^i$) & MPI\_Send($\bar {\vec u}^i$) \\
      % \bottomrule
      \noalign{\smallskip}\hline
    \end{tabular}
  \end{table}
% JBS Comment: Could you comment on why my_Step_Adjoint is described as
% additive, i.e. it uses +=?
% SG: Well I expected some comment like this... However, it's hard explain without going to deep into details for AD. I added a short explanation in the last paragraph and refer to the standard AD book for details. However, I'm not too happy with this explanation and will think about a better way to formulate it. Maybe I need to change the AD introduction a bit to make it more clear... I'll talk to an AD expert in my group next week and change it maybe.
% JBS: OK, I'll leave that up to you.

\subsection{Adjoint XBraid interface}
Only two of the user-defined routines $h_l$ are nonlinear, namely $\mathtt{my\_Step}$ and $\mathtt{my\_Access}$, such that their differentiated versions $\bar h_l$ are non-constant and contain problem specific functions. The adjoint XBraid interface therefore consists only of two additional routines that provide their derivatives:
  \begin{itemize}
    \item $\mathtt{my\_Step\_adjoint}$: This function corresponds to taking one adjoint time step backwards in time. Given a design $\rho$, a state variable ${\vec u}^{i-1}$ at time step $i-1$ and an adjoint input variable $\bar {\vec u}^i$ at time step $i$, it updates the adjoint variable $\bar {\vec u}^{i-1}$ at that time step $i-1$ and the gradient $\bar \rho$ according to
    \begin{align}
      \bar {\vec u}^{i-1} &\mathrel{+}=  \left( \partial_{{\vec u}^{i-1}} \Phi^i({\vec u}^{i-1}, \rho) \right)^T \bar {\vec u}^i \\
      \bar{\rho} &\mathrel{+}= \left( \partial_{\rho} \Phi^i({\vec u}^{i-1}, \rho) \right)^T \bar {\vec u}^i.
    \end{align}
    \item $\mathtt{my\_Access\_adjoint}$ updates the adjoint variable $\bar {\vec u}^i$ at {time step $i$} and the gradient $\bar{\rho}$ according to the partial derivatives of the instantaneous quantity $f(\vec u^i, \rho)$:
    \begin{align}
      \bar {\vec u}^i & \mathrel{+}=  \nabla_{{\vec u}^i} f({\vec u}^i, \rho) \\
      \bar{\rho} &\mathrel{+}= \nabla_{\rho} f({\vec u}^i, \rho).
    \end{align}
  \end{itemize}
Similar to the primal XBraid solver, the adjoint user interface is non-intrusive to existing adjoint time marching codes as they propagate the same sensitivities of the time stepper $\Phi^i$ and $f$ backwards through the time domain in a step-by-step manner.
Even though the adjoint XBraid solver has been derived utilizing techniques from AD, the adjoint user interface is not restricted to the use of AD for generating the desired derivatives of $\Phi$ and $f$. The adjoint interface rather enables integration of any standard unsteady adjoint solver into a parallel-in-time framework.
% To prepare an existing unsteady adjoint code for the adjoint XBraid solver, only the backwards time stepper and the derivative of the instantaneous quantity $f$ need to be wrapped into the adjoint routines $\mathtt{my\_Step\_adjoint}$ and $\mathtt{my\_Access\_adjoint}$.

\section{Numerical Results}\label{sec:Numerics}

The adjoint XBraid solver is validated on a model problem that mimics unsteady flow behind bluff bodies at low Reynolds numbers. Such flows typically contain two dominant flow regimes: In the near wake, periodic vortices are forming asymmetrically which shed into the far wake where they slowly dissipate, known as the Karman Vortex street. The Van-der-Pol oscillator, a nonlinear limit cycle ODE, is used to model the near wake oscillations while the far wake is modeled by an advection-diffusion equation whose upstream boundary condition is determined by the oscillations~\cite{kanamaru2007van}:
\begin{align}\label{Numerics_AdvectionDiffusion}
  \partial_t v(t,x) + a\partial_x v(t,x) - \mu \partial_{xx} v(t,x) &= 0 \qquad \forall x \in (0,1), t\in (0,T)   \\
  v(t,0) - \mu \partial_x v(t,0) &= z(t) \quad  \forall t\in(0,T) \\
  \partial_{xx}v(t,1) &= 0  \qquad  \forall t\in (0,T)  \\
  v(0,x) &= 1 \qquad \forall x \in [0,1],
\end{align}
while the advection term dominates with $a=1$ and a small diffusion parameter $\mu=10^{-5}$ is chosen. The Van-der-Pol oscillator determines the upstream boundary $z(t)$ with
\begin{align}\label{PinT:Numerics_VDP}
  \begin{pmatrix}
    \dot z(t) \\ \dot w(t)
  \end{pmatrix}  &=
  \begin{pmatrix}
   w(t) \\ -z(t) + \rho\left(1 - z(t)^2\right)w(t)
 \end{pmatrix}
 \quad \forall t \in (0, T),
\end{align}
using the initial condition $z(0) = w(0) = 1$. The parameter $\rho>0$ influences the nonlinear damping term and serves as design variable that determines the state $u(t,x)\coloneqq\left(z(t), w(t), v(t,x)\right)^T$ through the model equations. The objective function measures the space-time average of the solution
\begin{align}\label{PinT:Numerics_objectivefunction}
    J = \frac 1 T \int_0^T \|u(t,\cdot)\|_2^2 \, \mathrm{d}t.
\end{align}

The time domain is discretized into $N$ time steps with $t_i = i\Dt$ for $i=1, \dots, N$ while a fixed the final time $t_N = 30$ and varying $N$ and $\Dt$ are used to set up different problem sizes with $N=60000, 120000, 240000$ and corresponding $\Dt = 0.0005, 0.00025, 0.000125$, respectively.
The implicit Crank-Nicolson time marching scheme is chosen to approximate the transient term of the model equations while the spatial derivatives are approximated with a second order linear upwind scheme for the advection and central finite differencing for the diffusive term on the spatial grid $x_j = j\delta x, \, j=1, \dots,n$ with $\delta x = 0.01, n=100$.
The nonlinear equations at each time step are solved by applying functional iterations for $\vec u^i \coloneqq (z^i, w^i, v^i_1, \dots, v^i_n)\in \R^{2+n}$.
These iterations are wrapped into the core user routine $\mathtt{my\_Step}$ that moves a state $\vec u^{i-1}$ to the next time step.
The user routine $\mathtt{my\_Access}$ then evaluates the instantaneous quantity $\|\vec u^i\|^2_2$ at time step {$i$}.
The sensitivities $\left(\partial_{\vec u^{i-1},\rho} \Phi^i(\vec u^{i-1}, \rho)\right)^T\bar {\vec u}^i$ and $\nabla_{\vec u^i, \rho}f(\vec u^i, \rho)$ that are needed in the adjoint interface routines $\mathtt{my\_Step\_adjoint}$ and $\mathtt{my\_Access\_adjoint}$ are generated using the AD-Software CoDiPack for differentiating through their corresponding primal routines in reverse mode~\cite{CoDiPack}.
The time grid hierarchy for the multigrid iterations in the primal and adjoint XBraid solver use a coarsening factor of $m=4$ with a maximum of three time grid levels.
Adding more levels generates coarse grid time step sizes incompatible with the nonlinear time step solver, leading it to diverge.
However, even with three levels, a reasonable speedup can be demonstrated.  Implementing a more stable nonlinear time step routine, that would allow for more coarsening in time, is future work.
 The consideration of spatial coarsening is also future work, which would control the $\Dt / \delta x$ ratio on coarser grids, which often has the effect of making the nonlinear time step solver more stable.

% \begin{figure}
%   \includegraphics[width=.9\textwidth]{figures/VDP_N60000_phasespace.pdf}
%   \caption{Phase-space solution of the Van-der-Pol oscillator for increasing $\rho$.}
%   \label{fig:VDP_phasespace}
% \end{figure}

\begin{figure}
  \includegraphics[width=.9\textwidth]{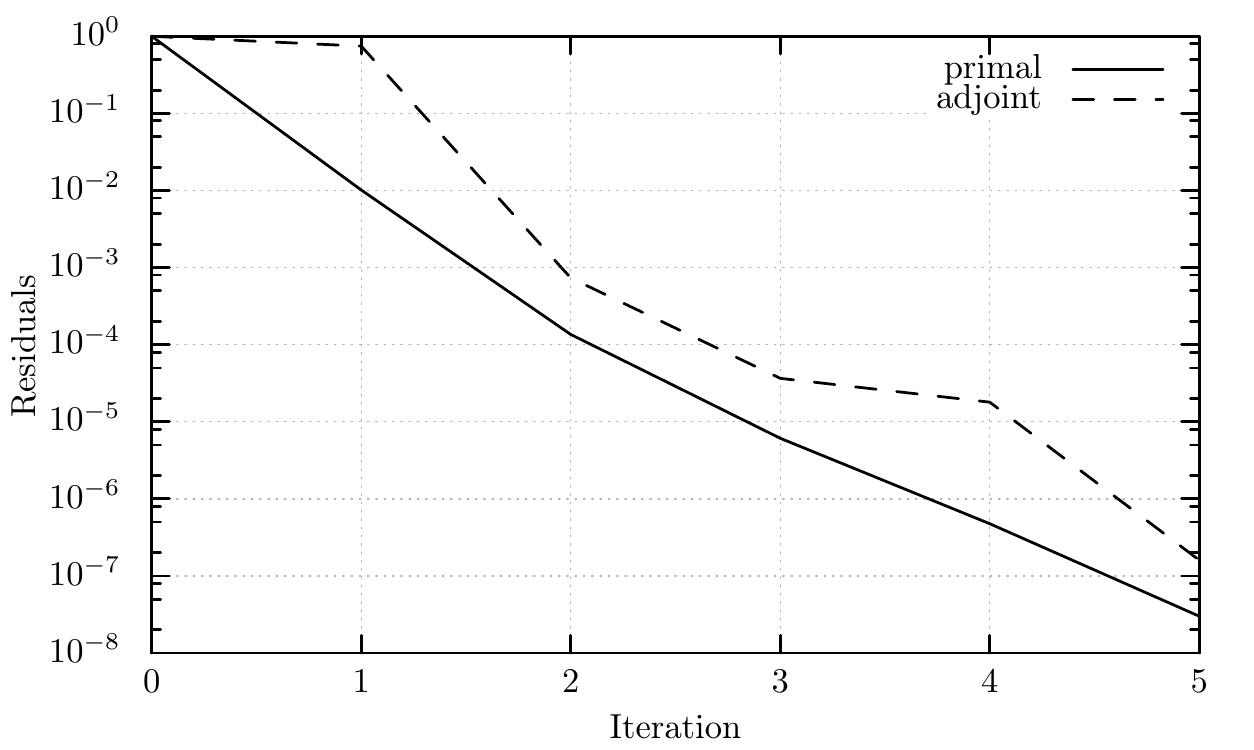}
  \caption{History of simultaneous piggyback iteration of primal and adjoint XBraid ($N=60000$).}
  \label{fig:piggybackiter}
\end{figure}

% \subsection{Piggyback iteration}
Figure \ref{fig:piggybackiter} shows the residuals of the primal and adjoint iterates, $\|\vec u_k - H(\vec u_k, \rho)\|_2$ and $\| \bar{\vec u}_k - \nabla_{\vec u}J_k - (\partial_{\vec u}H_k)^T\bar{\vec u}_k \|_2$ during a piggyback iteration as in \eqref{piggyback_itertion_primal}--\eqref{piggyback_itertion_adjoint} using a fixed design $\rho = 2$. As expected, both residuals drop simultaneously while the adjoint iterates exhibit a certain time lag.
The gradient provided by the adjoint XBraid solver is validated in Table \ref{tab:FinDiff_vs_XBraidAD} which shows good agreement with those computed from finite differences with a relative error below two percent.

\begin{table}
  \caption{Validation of the gradient computed from the adjoint XBraid solver with finite differences.}
  \label{tab:FinDiff_vs_XBraidAD}
% \caption{Table example}\label{tab:1}
\begin{tabular}{lllllll}
\hline\noalign{\smallskip}
        & $\rho$ & $\epsilon$ & finite differences & adjoint sensitivity & rel. error \\
\noalign{\smallskip}\hline\noalign{\smallskip}
$N=60000$ & $2$ & $10^{-6}$ & $0.231079186674421$ & $0.230724810643109$ & $0.15 \%$\\
$N=60000$ & $3$ & $10^{-6}$ & $0.224847819918494$ & $0.223017099689057$ & $0.81\%$\\
$N=60000$ & $4$ & $10^{-6}$ & $0.200429607133401$ & $0.196451436058937$ & $1.98\%$\\
\noalign{\smallskip}\hline\noalign{\smallskip}

$N=120000$ & $2$ & $10^{-6}$ & $0.231679529338891$ & $0.232395575709547$ & $0.31\%$\\
$N=120000$ & $3$ & $10^{-6}$  & $0.225517896579319$ & $0.225073911851428$ & $0.19\%$\\
$N=120000$ & $4$ & $10^{-6}$ & $0.199769404352068$ & $0.198695760225926$ & $0.54\%$\\
\noalign{\smallskip}\hline\noalign{\smallskip}

$N=240000$ & $2$ & $10^{-8}$ & $0.232875008165934$ & $0.234781281565824$ & $0.81\%$\\
$N=240000$ & $3$ & $10^{-8}$ & $0.227025065413500$ & $0.228025754956746$ & $0.44\%$\\
$N=240000$ & $4$ & $10^{-8}$ & $0.201409511291217$ & $0.204046297945359$ & $1.31 \%$\\

\noalign{\smallskip}\hline
\end{tabular}
\end{table}

\subsection{Parallel scaling for the primal and adjoint XBraid solver}

A weak scaling study for the primal and the adjoint XBraid solver is shown in Table \ref{tab:Primal_Adjoint_Iterations}. The reported speedups are computed by dividing the runtime of the primal and adjoint XBraid solver by that of a classical time-serial primal forward and adjoint backward time stepper, respectively.

% {Even though the adjoint runtimes show an overhead factor of $10$ when compared to the primal ones,}

{For these test cases, the runtimes of the adjoint XBraid solver show an overhead factor of about $12$ when compared to the corresponding runtimes of the primal one. However, a very similar factor of about $10$ can be observed for the time-serial computations. Thus, the observed overhead factor is mostly dominated by the adjoint time-stepping routine that computes $\left( \partial_{\vec u^{i-1},\rho}\Phi^i(\vec u^{i-1},\rho) \right)^T\vec{\bar u^i}$, which is needed in both the time-serial as well as the time-parallel computations. The additional overhead created by the AD-based derivation of the adjoint XBraid solver itself can therefore be estimated from $12/10 = 1.2$.}
% the fact that AD has been used to set up

\begin{table}
  \caption{Weak scaling study for the primal and adjoint XBraid solver.}
  \label{tab:Primal_Adjoint_Iterations}
  \begin{tabular}{ llllllr }
    \hline\noalign{\smallskip}
     &  & \phantom{a} & \multicolumn{4}{c}{Primal}  \\
    \cmidrule{4-7}
    Number of & Number of && Iterations & Runtime & Runtime  & Speedup  \\
    time steps & Cores && XBraid & XBraid & Serial  &  \\
    \noalign{\smallskip}\hline\noalign{\smallskip}
    60000  & 64   && 5  & 0.49 sec & 1.13 sec & 2.84 \\
    120000 & 128  && 4  & 0.50 sec & 2.97 sec & 5.97 \\
    240000 & 256  && 4  & 0.76 sec & 4.81 sec & 6.31 \\
    \noalign{\smallskip}\hline\noalign{\smallskip}
    &  & \phantom{a} & \multicolumn{4}{c}{Adjoint}  \\
   \cmidrule{4-7}
   Number of & Number of && Iterations & Runtime & Runtime  & Speedup  \\
   time steps & Cores && XBraid & XBraid & Serial  &  \\
    \noalign{\smallskip}\hline\noalign{\smallskip}
    60000  & 64   && 6  & 6.54 sec & 14.49 sec & 2.22 \\
    120000 & 128  && 4  & 6.12 sec & 26.27 sec & 4.23 \\
    240000 & 256  && 4  & 9.44 sec & 49.18 sec & 5.21 \\
    \noalign{\smallskip}\hline
  \end{tabular}
\end{table}

\begin{figure}
  \center
  \includegraphics[width=\textwidth]{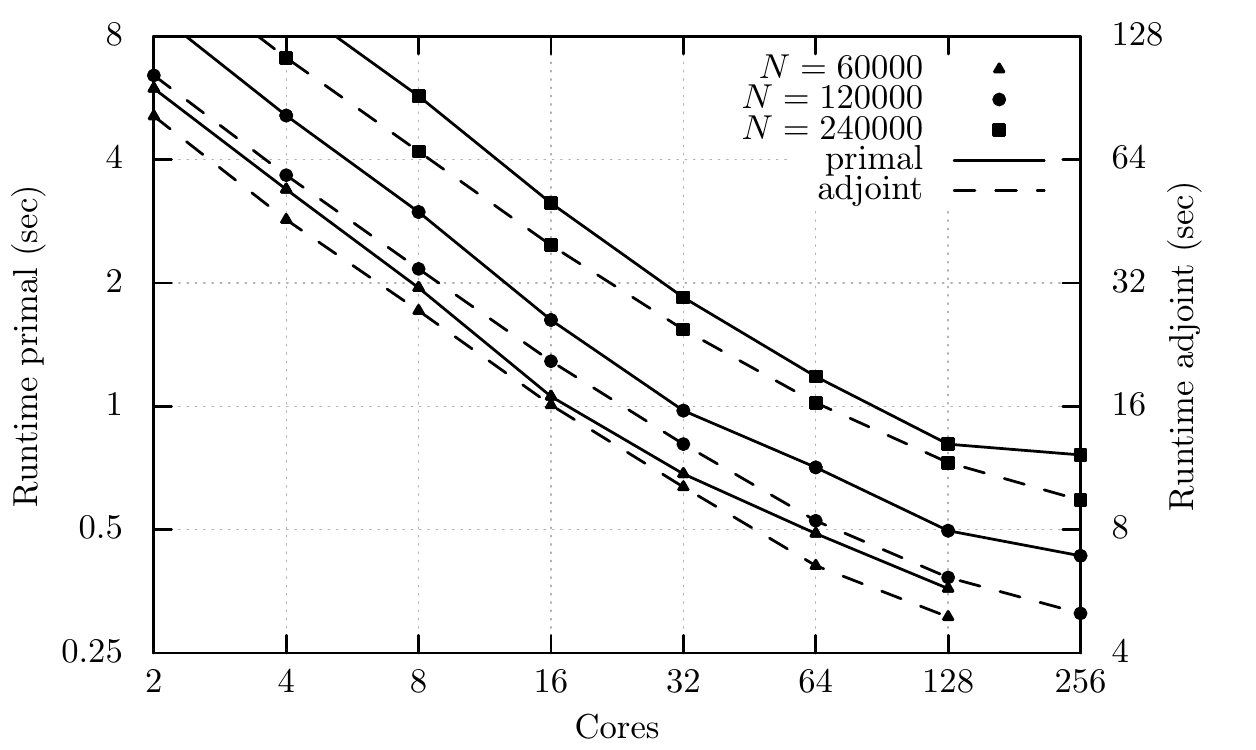}
  \caption{Strong scaling of primal (solid lines) and adjoint XBraid (dashed lines) solver for three different problem sizes.}
  \label{fig:PinT_primaladjoint_scaling}
\end{figure}

\begin{figure}
  \center
  \includegraphics[width=\textwidth]{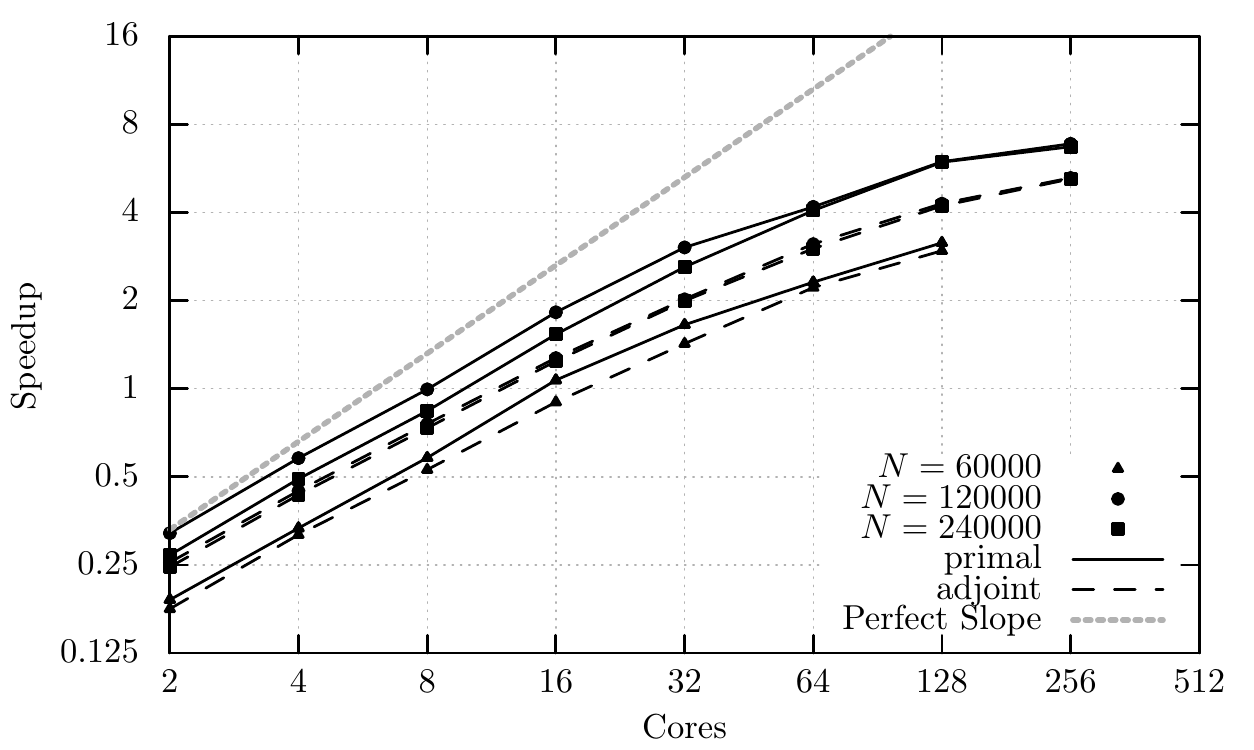}
  \caption{Speedup of primal and adjoint XBraid solver over serial forward and backward time stepping.}
  \label{fig:PinT_primaladjoint_speedup}
\end{figure}

Strong scaling results for primal as well as the adjoint runtimes are visualized in Figure \ref{fig:PinT_primaladjoint_scaling}{, where} the slope of reduction for the adjoint runtimes closely follows that of the primal one.
 This confirms that the adjoint XBraid solver indeed inherits the scaling behavior of the primal solver which is expected from its AD-based derivation. Since the primal XBraid solver is under active development, this property is particularly beneficial as improvements on primal scalability will automatically carry over to the adjoint code.
 The same data is used in Figure \ref{fig:PinT_primaladjoint_speedup} to plot speedup when compared to a time-serial primal and adjoint time marching scheme. It shows the potential of the time-parallel adjoint for speeding up the runtime of existing unsteady adjoint time marching schemes for sensitivity evaluation.


\section{Conclusion}

In this paper, we developed an adjoint solver that provides sensitivity computation for the parallel-in-time solver XBraid. XBraid applies nonlinear multigrid iterations to the time domain of unsteady partial differential equations and solves the resulting space-time system parallel-in-time. It operates through a high-level user interface that is non-intrusive to existing serial time marching schemes for solving unsteady PDEs. This property is of particular interest for applications where unsteady simulation tools have already been developed and refined over years, as is often the case for many CFD applications.

While the primal XBraid solver computes a space-time solution of the PDE for given input parameters and also evaluates objective functions that are of interest to the user, the proposed time-parallel adjoint XBraid solver computes derivatives of the objective function with respect to changes in the input parameters. Classical adjoint sensitivity computations for unsteady PDEs involve a forward-in-time sweep to compute the unsteady solution followed by a backwards-in-time loop to collect sensitivities. The parallel-in-time adjoint XBraid solver offers speedup by distributing the backwards-in-time phase onto multiple processors along the time domain. Its implementation is based on techniques of the reverse-mode of AD applied to one primal XBraid iteration. This yields a consistent discrete adjoint code that inherits parallel scaling properties from the primal solver and is non-intrusive to existing adjoint sequential time marching schemes.

The resulting adjoint solver adds two additional user routines to the primal XBraid interface: one for propagating sensitivities of the forward time stepper backwards in time and one for evaluating partial derivatives of the objective function at each time step. In cases where a time-serial unsteady adjoint solver is already available, this backwards time stepping capability can be easily wrapped into the adjoint XBraid interface with little extra coding.

The adjoint solver has been validated and tested on a model problem for advection-dominated flow. The original primal and adjoint time marching codes were limited to one processor such that a linear increase in the number of time steps results in a linear increase in corresponding runtime. This creates a situation analogous to the one where a spatially parallel code has reached its strong scaling limit. The parallel-in-time primal and adjoint XBraid solvers were able to achieve speedups of about $6$ (primal) and $5$ (adjoint) when compared to the serial ones, running on up to $256$ processors for the time parallelization. More importantly, {the scaling behavior of the adjoint code in this test case} is similar to that of the primal one such that improvements on the primal XBraid solver carry over to the adjoint implementation.

% \begin{table}
% % table caption is above the table
% \caption{Table example}\label{tab:1}
% % For LaTeX tables use
% \begin{tabular}{lll}
% \hline\noalign{\smallskip}
% first & second & third  \\
% \noalign{\smallskip}\hline\noalign{\smallskip}
% number & number & number \\
% number & number & number \\
% \noalign{\smallskip}\hline
% \end{tabular}
% \end{table}

\begin{acknowledgements}
The authors thanks Max Sagebaum (SciComp, TU Kaiserslautern) and Johannes Lotz (STCE, RWTH Aachen University) who provided insight and expertise on the implementation of AD.
\end{acknowledgements}

\bibliographystyle{spmpsci}      % mathematics and physical sciences
\bibliography{mybib}   % name your BibTeX data base

%\begin{thebibliography}{}
% etc
%\end{thebibliography}

\end{document}